\documentclass[a4paper,9pt]{article}
\usepackage{amsmath, amsfonts, amssymb}
\usepackage[english, russian]{babel}
\numberwithin{equation}{section}
\begin{document}
\sloppy

\bigskip

\begin{center}
\textbf{Decomposition formulas associated with the multivariable
\\
confluent hypergeometric functions}\\

Ergashev Ò.G.\\
\medskip
{ Institute of Mathematics, Uzbek
Academy of Sciences,  Tashkent, Uzbekistan. \\

{\verb ergashev.tukhtasin@gmail.com }}\\
\end{center}

The main object of this work is to show how some rather elementary
techniques based upon certain inverse pairs of symbolic operators would lead
us easily to several decomposition formulas associated with confluent
hypergeometric functions of two and more variables. Many operator identities
involving these pairs of symbolic operators are first constructed for this
purpose. By means of these operator identities several decomposition
formulas are found, which express the aforementioned hypergeometric
functions in terms of such simpler functions as the products of the Gauss
hypergeometric functions.


\section{Introduction}

A great interest in the theory of multiple hypergeometric
functions is motivated essentially by the fact that the solutions
of many applied problems involving, for example, partial
differential equations are obtainable with the help of such
hypergeometric functions (see, for details, \cite{SK}, p.47 et
seq., Section 1.7; see also the works \cite{{OS},{PS}} and the
references cited therein). For instance, the energy absorbed by
some nonferromagnetic conductor sphere included in an internal
magnetic field can be calculated with the help of such functions
\cite{L}. Hypergeometric functions in several variables are used
in physical and quantum chemical applications as well (cf.
\cite{{N},S}). Especially, many problems in gas dynamics lead to
solutions of degenerate second-order partial differential
equations which are then solvable in terms of multiple
hypergeometric functions.

We note that Riemann's functions and the fundamental solutions of
the degenerate second-order partial differential equations are
expressible by means of hypergeometric functions in several
variables \cite {{K},{SH1},{SH2},{SHC}}. In investigation of the
boundary-value problems for these partial differential equations,
we need decompositions for hypergeometric functions in several
variables in terms of simpler hypergeometric functions of the
Gauss and Appell types.

In addition to the Gaussian functions, which have received the
greatest attention in the literature, confluent functions have
been considered. For example, twenty confluent hypergeometric
functions of two variables exist; seven were introduced by Humbert
\cite{HM}, and the remaining ones by Horn \cite{HR} and by
Borng\"{a}sser \cite{B}. Certain confluent functions in three
variables were considered by Jain \cite{J} and by Exton \cite{EX},
but the entire set has not been given, i.e. confluent functions in
all directions of research have been little studied with respect
to other hypergeometric functions. While a brief account of such
functions is presented in \cite{SK}, we shall include the
definition of the important special class of confluent functions
and find the decomposition formulas for these functions.

Burchnall and Chaundy introduced the symbolic operators $\nabla $
and $\Delta $(see \cite{{BC1},{BC2}}) by means of which they
presented a number of expansion and decomposition formulas for
some double hypergeometric functions (only seven of which are
confluent functions) in terms of the classical Gauss
hypergeometric function of one variable. Recently Hasanov and
Srivastava \cite{{HS1},{HS2}} generalized the Burchnall-Chaundy's
operators and by making use of some technique based upon certain
inverse pairs of symbolic operators, the authors investigate
several decomposition formulas associated with Lauricella's (but
no confluent) hypergeometric functions of many variables when a
number of variables exceeds two.

In this paper we introduce other multivariable analog of Burchnall-Chaundy's
operators and find the decomposition formulas for some confluent
hypergeometric functions of two and more variables.

\section{Symbolic operators}

Burchnall and Chaundy \cite{{BC1},{BC2}}, and Chaundy \cite{C},
give a number of expansions of double hypergeometric functions in
series of simpler hypergeometric functions. Their method is based
upon the inverse pair of symbolic operators

\begin{equation}
\label{eq21} \nabla \left( {h} \right): = {\frac{{\Gamma \left(
{h} \right)\Gamma \left( {\delta _{1} + \delta _{2} + h}
\right)}}{{\Gamma \left( {\delta _{1} + h} \right)\Gamma \left(
{\delta _{2} + h} \right)}}}, \quad \Delta \left( {h} \right): =
{\frac{{\Gamma \left( {\delta _{1} + h} \right)\Gamma \left(
{\delta _{2} + h} \right)}}{{\Gamma \left( {h} \right)\Gamma
\left( {\delta _{1} + \delta _{2} + h} \right)}}},
\end{equation}
where
\begin{equation}
\label{eq22}
 \delta _{1} : = x_{1} {\frac{{\partial}} {{\partial x_{1}}} },\delta
_{2} : = x_{2} {\frac{{\partial}} {{\partial x_{2}}} }.
\end{equation}

The symbolic operators defined by (\ref{eq21}) and (\ref{eq22})
are limited only to functions of two variables,  therefore
recently Hasanov and Srivastava \cite {{HS1},{HS2}} generalized
these operators in the forms

\begin{equation}
\label{eq23} \tilde {\nabla} _{x_{1} :x_{2} ,...,x_{m}}  \left(
{h} \right): = {\frac{{\Gamma \left( {h} \right)\Gamma \left(
{\delta _{1} + ... + \delta _{m} + h} \right)}}{{\Gamma \left(
{\delta _{1} + h} \right)\Gamma \left( {\delta _{2} + ... + \delta
_{m} + h} \right)}}}
\end{equation}
and
\begin{equation}
\label{eq24} \tilde {\Delta} _{x_{1} :x_{2} ,...,x_{m}}  \left(
{h} \right): = {\frac{{\Gamma \left( {\delta _{1} + h}
\right)\Gamma \left( {\delta _{2} + ... + \delta _{m} + h}
\right)}}{{\Gamma \left( {h} \right)\Gamma \left( {\delta _{1} +
... + \delta _{m} + h} \right)}}},
\end{equation}
where

\[
\delta _{i} : = x_{i} {\frac{{\partial}} {{\partial x_{i}}}
}\,\,\,\left( {i = 1,...,m} \right).
\]

With the help of symbolic operators defined by (\ref{eq23}) and
(\ref{eq24}), decomposition formulas for many multiple
hypergeometric functions have been found. For example
\cite{{HS1},{HS2}},

\begin{equation}
\label{eq241}
\begin{array}{l}
  F_{\textbf{A}}^{(m)}\left(a,b_1,...,b_m;c_1,...,c_m;x_1,...,x_m\right) \\
= \tilde {\nabla} _{x_{1} :x_{2} ,...,x_{m}}  \left( {a}
\right)F\left(a,b_1;c_1;x_1\right)F_{\textbf{A}}^{(m-1)}\left(a,b_2,...,b_m;c_2,...,c_m;x_2,...,x_m\right) , \\
 \end{array}
\end{equation}

$$
\begin{array}{l}
  F_{\textbf{B}}^{(m)}\left(c_1,...,c_m,b_1,...,b_m;a; x_1,...,x_m\right) \\
= \tilde {\Delta} _{x_{1} :x_{2} ,...,x_{m}}  \left( {a}
\right)F\left(c_1,b_1;a;x_1\right)F_{\textbf{B}}^{(m-1)}\left(c_2,...,c_m,b_2,...,b_m;a;x_2,...,x_m\right) , \\
 \end{array}
$$
where
\begin{equation}
\label{eq31} F\left( {a,b,c;x} \right) = {\sum\limits_{m =
0}^{\infty}  {} }{\frac{{\left( {a} \right)_{m} \left( {b}
\right)_{m}}} {{\left( {c} \right)_{m}}} }{\frac{{x^{m}}}{{m!}}},
|x|<1,
\end{equation}
\begin{equation}
\label{eq41}
\begin{array}{l}
 F_{{\textbf{A}}}^{(m)} \left( {a,b_1,...,b_m;c_1,...,c_m;x_1,...,x_m} \right) \\=
{\sum\limits_{i_{1} ,...,i_{m} = 0}^{\infty}  {}} {\frac{{\left(
{a} \right)_{i_{1} + ... + i_{m}}  \left( {b_{1}}  \right)_{i_{1}}
...\left( {b_{m}}  \right)_{i_{m}}} } {{\left( {c_{1}}
\right)_{i_{1}}  ...\left( {c_{m}}  \right)_{i_{m}}
}}}{\frac{{x_{1}^{i_{1}}  ...x_{m}^{i_{m}}} } {{i_{1} !...i_{m}
!}}}, \,{\left| {x_{1}}  \right|} + {\left| {x_{2}}  \right|} +
... + {\left| {x_{m}}  \right|} < 1; \end{array}
\end{equation}

$$
\begin{array}{l}
 F_{{\textbf{B}}}^{(m)} \left( {c_1,...,c_m,b_1,...,b_m;a;x_1,...,x_m} \right)
 =\\
{\sum\limits_{i_{1} ,...,i_{m} = 0}^{\infty}  {}} {\frac{{\left(
{c_{1}}  \right)_{i_{1}}  ...\left( {c_{m}} \right)_{i_{m}} \left(
{b_{1}}  \right)_{i_{1}}  ...\left( {b_{m}} \right)_{i_{m}}} }
{{\left( {a} \right)_{i_{1} + ... + i_{m}}
}}}{\frac{{x_{1}^{i_{1}}  ...x_{m}^{i_{m}}} } {{i_{1} !...i_{m}
!}}}, \max \left( {\,{\left| {x_{1}}  \right|},{\left| {x_{2}}
\right|},...,{\left| {x_{m}}  \right|}} \right) < 1. \end{array}
$$

Here $ \left( {\mu } \right)_{k} : = \Gamma \left( {\mu + k}
\right) / \Gamma \left( {\mu} \right)$ is a Pochhammer symbol; $F
$ is Gauss hypergeometric function of one variable \cite[Chapter
2]{EM}; $F_{{\textbf{A}}}^{(m)}$ and $F_{{\textbf{B}}}^{(m)}$ are
multiple Lauricella hypergeometric functions \cite[p.115]{A}.

However, the recurrence of this formula did not allow further
advancement in the direction of increasing the number of
variables.

Further study of the properties of the Lauricella function defined
by (\ref{eq41}) showed that the formula (\ref{eq241}) can be
reduced to a more convenient form.

{\textbf{Lemma 1}\cite{ER}.\,\,\textit{The following formula holds
true at $m\in N$
\begin{equation}
\label{e28}\begin{array}{l} {F_{\textbf{A}}^{(m)} {\left(a,b_{1}
,....,b_{m} ;c_{1} ,....,c_{m} ;x_{1} ,...,x_{m} \right)}
 = {\sum\limits_{{\mathop {n_{i,j} = 0}\limits_{(2 \le i \le j \le m)}
}}^{\infty}  {{\frac{{(a)_{N_{2} (m,m)}}} {{{\mathop {n_{2,2}
!n_{2,3} ! \cdot \cdot \cdot n_{i,j} ! \cdot \cdot \cdot n_{m,m}
!}\limits_{(2 \le i \le j \le m)}}} }}}}} \hfill\\
 \cdot{\prod\limits_{k = 1}^{m} {{ {{\frac{{(b_{k} )_{M_{2} (k,m)}
}}{{(c_{k} )_{M_{2} (k,m)}}} x_{k}^{M_{2} (k,m)} F\left[a + N_{2}
(k,m),b_{k} + M_{2} (k,m);c_{k} + M_{2} (k,m);x_{k} \right]} }
}}},
\end{array}
\end{equation}
where
$$
M_{l} (k,m) = {\sum\limits_{i = l}^{k} {n_{i,k} +}}
{\sum\limits_{i = k + 1}^{m} {n_{k + 1,i}}},\,\, \quad N_{l} (k,m)
= {\sum\limits_{i = l}^{k + 1} {{\sum\limits_{j = i}^{m}
{n_{i,j}}} } } , \, l \in N.
$$
}

It should be noted that the symbolic operators $\delta _{1} $ and
$\delta _{2} $ defined by (\ref{eq22}) in the one-dimensional case
take the form $\delta : = xd / dx$ and such an operator is used in
solving problems of the operational calculus \cite[p.26]{P}.

We now introduce here the other multivariable analogues of the
Burchnall-Chaundy symbolic operators $\nabla \left( {h} \right)$
and $\Delta \left( {h} \right)$ defined by (\ref{eq21}):

\begin{equation}
\label{eq25} \tilde {\nabla} _{x;y}^{m,n} \left( {h} \right): =
{\frac{{\Gamma \left( {h} \right)\Gamma \left( {h + \delta _{1} +
... + \delta _{m} - \sigma _{1} - ... - \sigma _{n}}
\right)}}{{\Gamma \left( {h + \delta _{1} + ... + \delta _{m}}
\right)\Gamma \left( {h - \sigma _{1} - ... - \sigma _{n}}
\right)}}}
\end{equation}

\begin{equation}
\label{eq26}
 = {\sum\limits_{k = 0}^{\infty}  {{\frac{{\left( { - \delta _{1} - ... -
\delta _{m}}  \right)_{k} \left( {\sigma _{1} + ... + \sigma _{n}}
\right)_{k}}} {{(h)_{k} k!}}}}},
\end{equation}

\begin{equation}
\label{eq27} \tilde {\Delta} _{x;y}^{m,n} \left( {h} \right): =
{\frac{{\Gamma \left( {h + \delta _{1} + ... + \delta _{m}}
\right)\Gamma \left( {h - \sigma _{1} - ... - \sigma _{n}}
\right)}}{{\Gamma \left( {h} \right)\Gamma \left( {h + \delta _{1}
+ ... + \delta _{m} - \sigma _{1} - ... - \sigma _{n}} \right)}}}
\end{equation}

\begin{equation}
\label{eq28}
 = {\sum\limits_{k = 0}^{\infty}  {{\frac{{\left( {\delta _{1} + ... +
\delta _{m}}  \right)_{k} \left( { - \sigma _{1} - ... - \sigma
_{n}} \right)_{k}}} {{(1 - h)_{k} k!}}}}},
\end{equation}
where

\begin{equation}
\label{eq29} x: = \left( {x_{1} ,...,x_{m}}  \right), y: = \left(
{y_{1} ,...,y_{n}}  \right),
\end{equation}

$$
 \delta _{i} : = x_{i} {\frac{{\partial}} {{\partial
x_{i}}} }\,\,, \sigma _{j} : = y_{j} {\frac{{\partial}} {{\partial
y_{j}}} }, \quad i = 1,...,m,\,\,j = 1,...,n;\,\,m,n \in {\rm N}.
$$

In addition, we consider operators which are equal to the
Hasanov-Srivastava's symbolic operators $\tilde {\nabla} \left(
{h} \right)$ and $\tilde {\Delta} \left( {h} \right)$ defined by
(\ref{eq23}) and (\ref{eq24}):

$$ \tilde {\nabla} _{x; -} ^{m,0} \left( {h} \right): = \tilde
{\nabla} _{x_{1} :x_{2} ,...,x_{m}}  \left( {h} \right), \quad
\tilde {\Delta} _{x; -} ^{m,0} \left( {h} \right): = \tilde
{\Delta} _{x_{1} :x_{2} ,...,x_{m}}  \left( {h} \right), m \in
{\rm N};
$$

$$\tilde {\nabla} _{ - ;y}^{0,n} \left( {h} \right): = \tilde
{\nabla} _{ - y_{1} : - y_{2} ,..., - y_{n}}  \left( {h} \right),
\quad \tilde {\Delta} _{ - ;y}^{0,n} \left( {h} \right): = \tilde
{\Delta} _{ - y_{1} : - y_{2} ,..., - y_{n}}  \left( {h} \right),
n \in {\rm N}.
$$

It is obvious that

$$
\tilde {\nabla} _{x; -} ^{1,0} \left( {h} \right) = \tilde {\Delta} _{x; -
}^{1,0} \left( {h} \right) = \tilde {\nabla} _{ - ;y}^{0,1} \left( {h}
\right) = \tilde {\Delta} _{ - ;y}^{0,1} \left( {h} \right) = 1.
$$

\textbf{Lemma 2.} \,\,\textit{Let be $f: = f\left( {x,y} \right)$
function with variables $x$ and $y$ in (\ref{eq29}). Then
following equalities hold true for any $m,\,n \in {\rm N}$:}

\begin{equation}
\label{eq213} \left( { - {\sum\limits_{i = 1}^{m} {x_{i}
{\frac{{\partial}} {{\partial x_{i}}} }}}}  \right)_{k} f = ( -
1)^{k}k!{\sum\limits_K {{\prod\limits_{s = 1}^{m}
{{\frac{{x_{s}^{i_{s}}} } {{i_{s} !}}}}}  \cdot {\frac{{\partial
^{k}f}}{{\partial x_{1}^{i_{1}}  ...\partial x_{m}^{i_{m}}} } }}}
,\,\,\,k \in {\rm N} \cup {\left\{ {0} \right\}};
\end{equation}

\begin{equation}
\label{eq214} \left( {{\sum\limits_{j = 1}^{n} {y_{j}
{\frac{{\partial}} {{\partial y_{j} }}}}}}  \right)_{k} f =
{\left\{ {{\begin{array}{*{20}c}

{f,\,\,\,\,\,\,\,\,\,\,\,\,\,\,\,\,\,\,\,\,\,\,\,\,\,\,\,\,\,\,\,\,\,\,\,\,\,\,\,\,\,\,\,\,\,\,k
= 0,} \hfill \\
 {k!{\sum\limits_{l = 1}^{k} {\left( {{\begin{array}{*{20}c}
 {l - 1} \hfill \\
 {k - 1} \hfill \\
\end{array}}}  \right)}} {\sum\limits_L {{\prod\limits_{s = 1}^{n} {{\frac{{y_{s}^{j_{s}}} } {{j_{s}
!}}}}}  \cdot {\frac{{\partial ^{l}f}}{{\partial y_{1}^{j_{1}}  ...\partial
y_{n}^{j_{n}}} } }}} ,\,\,k \in {\rm N},} \hfill \\
\end{array}}}  \right.}
\end{equation}
\textit{where}
$$ K:=\left\{\left(i_1,...,i_m\right): i_1\geq0,...,i_m\geq0, i_1+...+i_m=k\right\},$$
 $$ L:=\left\{\left(j_1,...,j_n\right): j_1\geq0,...,j_n\geq0, j_1+...+j_n=l\right\}.$$

The lemma 2 is proved by method of mathematical induction.

\section{Decomposition formulas in the two-variable case}

In this section we shall give the decomposition formulas for the
following  hypergeometric functions of two variables
\cite[pp.225-226]{EM}:

$$
 H_{2} \left( {a,b,c,e,d;x,y} \right) =
{\sum\limits_{m,n = 0}^{\infty}  {} }{\frac{{\left( {a} \right)_{m
- n} \left( {b} \right)_{m} \left( {c} \right)_{n} \left( {e}
\right)_{n}}} {{\left( {d} \right)_{m} m!n!}}}x^{m}y^{n},
$$

\begin{equation}
\label{eq35} {\rm H}_{2} \left( {a,b,c;d;x,y} \right) =
{\sum\limits_{m,n = 0}^{\infty} {{\frac{{\left( {a} \right)_{m -
n} \left( {b} \right)_{m} \left( {c} \right)_{n}}} {{\left( {d}
\right)_{m} m!n!}}}x^{m}y^{n}}} ,\,\,{\left| {x} \right|} < 1,
\end{equation}

\begin{equation}
\label{eq36} {\rm H}_{3} \left( {a,b;d;x,y} \right) =
{\sum\limits_{m,n = 0}^{\infty} {{\frac{{\left( {a} \right)_{m -
n} \left( {b} \right)_{m}}} {{\left( {d} \right)_{m}
m!n!}}}x^{m}y^{n}}} ,\,\,{\left| {x} \right|} < 1,
\end{equation}

\begin{equation}
\label{eq37} {\rm H}_{4} \left( {a,b;d;x,y} \right) =
{\sum\limits_{m,n = 0}^{\infty} {{\frac{{\left( {a} \right)_{m -
n} \left( {b} \right)_{n}}} {{\left( {d} \right)_{m}
m!n!}}}x^{m}y^{n}}} ,
\end{equation}

\begin{equation}
\label{eq38} {\rm H}_{5} \left( {a;d;x,y} \right) =
{\sum\limits_{m,n = 0}^{\infty} {{\frac{{\left( {a} \right)_{m -
n}}} {{\left( {d} \right)_{m} m!n!}}}x^{m}y^{n}}} ,\,\,
\end{equation}

\begin{equation}
\label{eq39} {\rm H}_{11} \left( {a,b,c;d;x,y} \right) =
{\sum\limits_{m,n = 0}^{\infty} {{\frac{{\left( {a} \right)_{m -
n} \left( {b} \right)_{n} \left( {c} \right)_{n}}} {{\left( {d}
\right)_{m} m!n!}}}x^{m}y^{n}}} ,\,\,{\left| {y} \right|} < 1,
\end{equation}
where $a,b,c,d,e$ are complex numbers, $d \ne 0, - 1, - 2,...$. We
note that hypergeometric functions defined by
(\ref{eq35})-(\ref{eq39}) are confluent functions (In the
literature it is customary to denote the confluent functions
through the capital letters of the Greek alphabet).

In the special case when $m = n = 1,$ the symbolic operators
(\ref{eq25})-(\ref{eq28}) and equalities
(\ref{eq213})-(\ref{eq214}) take a simpler forms:

\begin{equation}
\label{eq310} \tilde {\nabla} _{x,y}^{1,1} \left( {h} \right): =
{\frac{{\Gamma \left( {h} \right)\Gamma \left( {h + \delta -
\sigma}  \right)}}{{\Gamma \left( {h + \delta}  \right)\Gamma
\left( {h - \sigma}  \right)}}}  = {\sum\limits_{k = 0}^{\infty}
{{\frac{{\left( { - \delta}  \right)_{k} \left( {\sigma}
\right)_{k}}} {{(h)_{k} k!}}}}} ,
\end{equation}

\begin{equation}
\label{eq311} \tilde {\Delta} _{x,y}^{1,1} \left( {h} \right): =
{\frac{{\Gamma \left( {h + \delta}  \right)\Gamma \left( {h -
\sigma}  \right)}}{{\Gamma \left( {h} \right)\Gamma \left( {h +
\delta - \sigma}  \right)}}} = {\sum\limits_{k = 0}^{\infty}
{{\frac{{\left( {\delta}  \right)_{k} \left( { - \sigma}
\right)_{k}}} {{(1 - h)_{k} k!}}}}} ,
\end{equation}

\begin{equation}
\label{eq312} \left( { - x{\frac{{\partial}} {{\partial x}}}}
\right)_{k} f = ( - 1)^{k}x^{k}{\frac{{\partial ^{k}f}}{{\partial
x^{k}}}},\,\,\,k \in {\rm N} \cup {\left\{ {0} \right\}},
\end{equation}

\begin{equation}
\label{eq313} \left( {y{\frac{{\partial}} {{\partial y}}}}
\right)_{k} f = {\left\{ {{\begin{array}{*{20}c}

{f,\,\,\,\,\,\,\,\,\,\,\,\,\,\,\,\,\,\,\,\,\,\,\,\,\,\,\,\,\,\,\,\,\,\,\,\,\,\,\,\,\,\,\,\,\,\,k
= 0,} \hfill \\
 {k!{\sum\limits_{l = 1}^{k} {\left( {{\begin{array}{*{20}c}
 {l - 1} \hfill \\
 {k - 1} \hfill \\
\end{array}}}  \right)}} {\frac{{y^{l}}}{{l!}}}{\frac{{\partial
^{l}f}}{{\partial y^{l}}}},\,\,k \in {\rm N},\,\,\,\,} \hfill \\
\end{array}}}  \right.}
\end{equation}
where
\[
\delta : = x{\frac{{\partial}} {{\partial x}}}\,\,,
\quad
\sigma : = y{\frac{{\partial}} {{\partial y}}}.
\]

By applying the pair of symbolic operators (\ref{eq310}) and
(\ref{eq311}), we find the following set of operator identities:

\begin{equation}
\label{eq314} H_{2} \left( {a,b,c,d;e;x,y} \right) = \tilde
{\nabla} _{x,y}^{1,1} \left( {a} \right)F\left( {a,b;e;x}
\right)F\left( {c,d;1 - a; - y} \right),
\end{equation}

\begin{equation}
\label{eq315} F\left( {a,b;e;x} \right)F\left( {c,d;1 - a;y}
\right) = \tilde {\Delta }_{x,y}^{1,1} \left( {a} \right)H_{2}
\left( {a,b,c,d;e;x, - y} \right),
\end{equation}

\begin{equation}
\label{eq316} {\rm H}_{2} \left( {a,b,c;d;x,y} \right) = \tilde
{\nabla} _{x,y}^{1,1} \left( {a} \right)F\left( {a,b;d;x}
\right){}_{1}F_{1} \left( {c;1 - a; - y} \right),
\end{equation}

\begin{equation}
\label{eq317} F\left( {a,b;d;x} \right){}_{1}F_{1} \left( {c;1 -
a; - y} \right) = \tilde {\Delta} _{x,y}^{1,1} \left( {a}
\right){\rm H}_{2} \left( {a,b,c;d;x,y} \right),
\end{equation}

\begin{equation}
\label{eq318} {\rm H}_{3} \left( {a,b;d;x,y} \right) = \tilde
{\nabla} _{x,y}^{1,1} \left( {a} \right)F\left( {a,b;d;x}
\right){}_{0}F_{1} \left( {1 - a; - y} \right),
\end{equation}

\begin{equation}
\label{eq319} F\left( {a,b;d;x} \right){}_{0}F_{1} \left( {1 -
a;y} \right) = \tilde {\Delta} _{x,y}^{1,1} \left( {a} \right){\rm
H}_{3} \left( {a,b;d;x, - y} \right),
\end{equation}

\begin{equation}
\label{eq320} {\rm H}_{4} \left( {a,b;d;x,y} \right) = \tilde
{\nabla} _{x,y}^{1,1} \left( {a} \right){}_{1}F_{1} \left( {a;d;x}
\right){}_{1}F_{1} \left( {b;1 - a; - y} \right),
\end{equation}

\begin{equation}
\label{eq321} {}_{1}F_{1} \left( {a;d;x} \right){}_{1}F_{1} \left(
{b;1 - a;y} \right) = \tilde {\Delta} _{x,y}^{1,1} \left( {a}
\right){\rm H}_{4} \left( {a,b;d;x, - y} \right),
\end{equation}

\begin{equation}
\label{eq322} {\rm H}_{5} \left( {a;d;x,y} \right) = \tilde
{\nabla} _{x,y}^{1,1} \left( {a} \right){}_{1}F_{1} \left( {a;d;x}
\right){}_{0}F_{1} \left( {1 - a; - y} \right),
\end{equation}

\begin{equation}
\label{eq323} {}_{1}F_{1} \left( {a;d;x} \right){}_{0}F_{1} \left(
{1 - a;y} \right) = \tilde {\Delta} _{x,y}^{1,1} \left( {a}
\right){\rm H}_{5} \left( {a;d;x, - y} \right),
\end{equation}

\begin{equation}
\label{eq324} {\rm H}_{11} \left( {a,b,c;d;x,y} \right) = \tilde
{\nabla} _{x,y}^{1,1} \left( {a} \right){}_{1}F_{1} \left( {a;d;x}
\right)F\left( {b,c;1 - a; - y} \right),
\end{equation}

\begin{equation}
\label{eq325} {}_{1}F_{1} \left( {a;d;x} \right)F\left( {b,c;1 -
a; - y} \right) = \tilde {\Delta} _{x,y}^{1,1} \left( {a}
\right){\rm H}_{11} \left( {a,b,c;d;x,y} \right),
\end{equation}
where
$${}_pF_q\left(a_1,...,a_p;b_1,...,b_q;x\right):=\sum\limits_{n=0}^\infty\frac{(a_1)_n...(a_p)_n}{(b_1)_n...(b_q)_n}\frac{x^n}{n!}, |x|<1,$$
is a generalized Gauss hypergeometric function, $F$ is the famous
Gauss function defined by (\ref{eq31}).

By using equalities (\ref{eq312}) and (\ref{eq313}) from the
operator identities (\ref{eq314}) to (\ref{eq325}) we can derive
the following decomposition formulas for double hypergeometric
functions $H_{2}$ ,  ${\rm H}_{2} - {\rm H}_{5} $ and ${\rm
H}_{11} :$

\begin{equation}
\label{eq326}
\begin{array}{l}
 H_{2} \left( {a,b,c,d;e;x,y} \right) = F\left( {a,b;e;x} \right)F\left(
{c,d;1 - a; - y} \right)\, + {\sum\limits_{k = 1}^{\infty}
{{\sum\limits_{l = 1}^{k} {}} {\frac{{\left( { - 1} \right)^{k +
l}\left( {k - 1}
\right)!}}{{\left( {l - 1} \right)!l!\left( {k - l} \right)!}}}}}  \\
 \cdot {\frac{{\left( {b} \right)_{k} \left( {c} \right)_{l} \left( {d}
\right)_{l}}} {{\left( {1 - a} \right)_{l} \left( {e} \right)_{k}
}}}x^{k}y^{l}F\left( {a + k,b + k;e + k;x} \right)F\left( {c + l,d
+ l;1 - a
+ l; - y} \right), \\
 \end{array}
\end{equation}

\begin{equation}
\label{eq327}
\begin{array}{l}
 F\left( {a,b;e;x} \right)F\left( {c,d;1 - a;y} \right) = H_{2} \left(
{a,b,c,d;e;x, - y} \right) + {\sum\limits_{k = 1}^{\infty}
{{\sum\limits_{l = 1}^{k} {{\frac{{\left( {k - 1}
\right)!}}{{\left( {l - 1} \right)!l!\left(
{k - l} \right)!}}}}}} }  \\
 \cdot {\frac{{\left( { - 1} \right)^{k - l}\left( {b} \right)_{l} \left(
{c} \right)_{k} \left( {d} \right)_{k}}} {{\left( {1 - a}
\right)_{k} \left( {1 - a} \right)_{k - l} \left( {e}
\right)_{l}}} }x^{l}y^{k}H_{2} \left( {a
- k + l,b + l,c + k,d + k;e + l;x, - y} \right), \\
 \end{array}
\end{equation}

\begin{equation}
\label{eq328}
\begin{array}{l}
 {\rm H}_{2} \left( {a,b,c;d;x,y} \right) = F\left( {a,b;d;x}
\right){}_{1}F_{1} \left( {c;1 - a; - y} \right) + {\sum\limits_{k
= 1}^{\infty}  {{\sum\limits_{l = 1}^{k} {{\frac{{( - 1)^{k +
l}\left( {k - 1}
\right)!}}{{\left( {l - 1} \right)!l!\left( {k - l} \right)!}}}}}} }  \\
 \cdot {\frac{{\left( {b} \right)_{k} \left( {c} \right)_{l}}} {{\left( {1 -
a} \right)_{l} \left( {d} \right)_{k}}} }x^{k}y^{l}F\left( {a +
k,b + k;d +
k;x} \right){}_{1}F_{1} \left( {c + l;1 - a + l; - y} \right), \\
 \end{array}
\end{equation}

\begin{equation}
\label{eq329}
\begin{array}{l}
 F\left( {a,b;d;x} \right){}_{1}F_{1} \left( {c;1 - a;y} \right) = {\rm
H}_{2} \left( {a,b,c;d;x, - y} \right) + {\sum\limits_{k =
1}^{\infty} {{\sum\limits_{l = 1}^{k} {{\frac{{\left( { - 1}
\right)^{k - l}\left( {k -
1} \right)!}}{{\left( {l - 1} \right)!l!\left( {k - l} \right)!}}}}}} }  \\
 \cdot {\frac{{\left( {b} \right)_{l} \left( {c} \right)_{k}}} {{\left( {1 -
a} \right)_{k} \left( {1 - a} \right)_{k - l} \left( {d}
\right)_{l} }}}x^{l}y^{k}{\rm H}_{2} \left( {a - k + l,b + l,c +
k;d + l;x, - y}
\right), \\
 \end{array}
\end{equation}

\begin{equation}
\label{eq330}
\begin{array}{l}
 {\rm H}_{3} \left( {a,b;d;x,y} \right) = F\left( {a,b;d;x}
\right){}_{0}F_{1} \left( {1 - a; - y} \right) + {\sum\limits_{k =
1}^{\infty}  {{\sum\limits_{l = 1}^{k} {{\frac{{( - 1)^{k +
l}\left( {k - 1}
\right)!}}{{(l - 1)!l!\left( {k - l} \right)!}}}}}} }  \\
 \cdot {\frac{{(b)_{k}}} {{(1 - a)_{l} (d)_{k}}} }x^{k}y^{l}F\left( {a + k,b
+ k;d + k;x} \right){}_{0}F_{1} \left( {1 - a + l; - y} \right), \\
 \end{array}
\end{equation}

\begin{equation}
\label{eq331}
\begin{array}{l}
 F\left( {a,b;d;x} \right){}_{0}F_{1} \left( {1 - a;y} \right) = {\rm H}_{3}
\left( {a,b;d;x, - y} \right) + {\sum\limits_{k = 1}^{\infty}
{{\sum\limits_{l = 1}^{k} {{\frac{{\left( { - 1} \right)^{k -
l}\left( {k -
1} \right)!}}{{\left( {l - 1} \right)!l!\left( {k - l} \right)!}}}}}} }  \\
 \cdot {\frac{{\left( {b} \right)_{l}}} {{\left( {1 - a} \right)_{k} \left(
{1 - a} \right)_{k - l} \left( {d} \right)_{l}}} }x^{l}y^{k}{\rm
H}_{3}
\left( {a - k + l,b + l;d + l;x, - y} \right), \\
 \end{array}
\end{equation}

\begin{equation}
\label{eq332}
\begin{array}{l}
 {\rm H}_{4} \left( {a,b;d;x,y} \right) = {}_{1}F_{1} \left( {a;d;x}
\right){}_{1}F_{1} \left( {b;1 - a; - y} \right) + {\sum\limits_{k
= 1}^{\infty}  {{\sum\limits_{l = 1}^{k} {{\frac{{( - 1)^{k +
l}\left( {k - 1}
\right)!}}{{\left( {l - 1} \right)!l!\left( {k - l} \right)!}}}}}} }  \\
 \cdot {\frac{{\left( {b} \right)_{l}}} {{\left( {1 - a} \right)_{l} \left(
{d} \right)_{k}}} }x^{k}y^{l}{}_{1}F_{1} \left( {a + k;d + k;x}
\right){}_{1}F_{1} \left( {b + l;1 - a + l; - y} \right), \\
 \end{array}
\end{equation}

\begin{equation}
\label{eq333}
\begin{array}{l}
 {}_{1}F_{1} \left( {a;d;x} \right){}_{1}F_{1} \left( {b;1 - a;y} \right) =
{\rm H}_{4} \left( {a,b;d;x, - y} \right) + {\sum\limits_{k =
1}^{\infty} {{\sum\limits_{l = 1}^{k} {}} {\frac{{\left( { - 1}
\right)^{k - l}\left( {k
- 1} \right)!}}{{\left( {l - 1} \right)!l!\left( {k - l} \right)!}}}}}  \\
 \cdot {\frac{{\left( {b} \right)_{k}}} {{\left( {1 - a} \right)_{k} \left(
{1 - a} \right)_{k - l} \left( {d} \right)_{l}}} }x^{l}y^{k}{\rm
H}_{4}
\left( {a - k + l,b + k;d + l;x, - y} \right), \\
 \end{array}
\end{equation}

\begin{equation}
\label{eq334}
\begin{array}{l}
 {\rm H}_{5} \left( {a;d;x,y} \right) = {}_{1}F_{1} \left( {a;d;x}
\right){}_{0}F_{1} \left( {1 - a; - y} \right) + {\sum\limits_{k =
1}^{\infty}  {{\sum\limits_{l = 1}^{k} {}} {\frac{{( - 1)^{k +
l}\left( {k -
1} \right)!}}{{\left( {l - 1} \right)!l!\left( {k - l} \right)!}}}}}  \\
 \cdot {\frac{{1}}{{\left( {1 - a} \right)_{l} \left( {d} \right)_{k}
}}}x^{k}y^{l}{}_{1}F_{1} \left( {a + k;d + k;x} \right){}_{0}F_{1}
\left( {1
- a + l; - y} \right), \\
 \end{array}
\end{equation}

\begin{equation}
\label{eq335}
\begin{array}{l}
 {}_{1}F_{1} \left( {a;d;x} \right){}_{0}F_{1} \left( {1 - a;y} \right) =
{\rm H}_{5} \left( {a,b;d;x, - y} \right) + {\sum\limits_{k =
1}^{\infty} {{\sum\limits_{l = 1}^{k} {}} {\frac{{\left( { - 1}
\right)^{k - l}\left( {k
- 1} \right)!}}{{\left( {l - 1} \right)!l!\left( {k - l} \right)!}}}}}  \\
 \cdot {\frac{{1}}{{\left( {1 - a} \right)_{k} \left( {1 - a} \right)_{k -
l} \left( {d} \right)_{l}}} }x^{l}y^{k}{\rm H}_{5} \left( {a - k +
l;d +
l;x, - y} \right), \\
 \end{array}
\end{equation}

\begin{equation}
\label{eq336}
\begin{array}{l}
 {\rm H}_{11} \left( {a,b,c;d;x,y} \right) = {}_{1}F_{1} \left( {a;d;x}
\right)F\left( {b,c;1 - a; - y} \right) + {\sum\limits_{k =
1}^{\infty} {{\sum\limits_{l = 1}^{k} {}} {\frac{{\left( { - 1}
\right)^{k + l}\left( {k
- 1} \right)!}}{{\left( {l - 1} \right)!l!\left( {k - l} \right)!}}}}}  \\
 \cdot {\frac{{\left( {b} \right)_{l} \left( {c} \right)_{l}}} {{\left( {1 -
a} \right)_{l} \left( {d} \right)_{k}}} }x^{k}y^{l}{}_{1}F_{1}
\left( {a +
k;d + k;x} \right)F\left( {b + l,c + l;1 - a + l; - y} \right), \\
 \end{array}
\end{equation}

\begin{equation}
\label{eq337}
\begin{array}{l}
 {}_{1}F_{1} \left( {a;d;x} \right)F\left( {b,c;1 - a;y} \right) = {\rm
H}_{11} \left( {a,b,c;d;x, - y} \right) + {\sum\limits_{k =
1}^{\infty} {{\sum\limits_{l = 1}^{k} {}} {\frac{{\left( { - 1}
\right)^{k - l}\left( {k
- 1} \right)!}}{{\left( {l - 1} \right)!l!\left( {k - l} \right)!}}}}}  \\
 \cdot {\frac{{\left( {b} \right)_{k} \left( {c} \right)_{k}}} {{\left( {1 -
a} \right)_{k} \left( {1 - a} \right)_{k - l} \left( {d}
\right)_{l} }}}x^{l}y^{k}{\rm H}_{11} \left( {a - k + l,b + k,c +
k;d + l;x, - y}
\right). \\
 \end{array}
\end{equation}

The expansions (\ref{eq326})-(\ref{eq337}) can be proved without
symbolic methods by comparing coefficients of equal powers of $x$
and $y$ on both sides.

\section{Decomposition formulas for the
multivariable confluent hypergeometric function }

An interesting unification (and generalization) of multiple
Lauricella's functions $\,F_{{\textbf{A}}}^{(m)} $ and
$F_{{\textbf{B}}}^{(m)} $ and Horn's functions of two variables
$H_{2} $ was considered by Erd\'{e}lyi \cite{E}(see also
\cite[p.74]{SK}), who defined his general function in the form:

\begin{equation}
\label{eq44}
\begin{array}{l}
 H_{m + n,m} \left( {a,b_1,...,b_m,d_{1} ,...,d_{n} ,e_{1} ,...,e_{n} ;c_1,...,c_m;x,y}
 \right)\\
= {\sum {} }
 {\frac{(a)_{i_{1} + ... + i_{m} - j_{1} - ... - j_{n}}{\left( {b_{1}}  \right)_{i_{1}}  ...\left( {b_{m}}
\right)_{i_{m}}  \left( {d_{1}}  \right)_{j_{1}}  ...\left(
{d_{n}} \right)_{j_{n}}  \left( {e_{1}}  \right)_{j_{1}} ...\left(
{e_{n}} \right)_{j_{n}}} } {{\left( {c_{1}} \right)_{i_{1}}
...\left( {c_{m}} \right)_{i_{m}}} } }{\frac{{x_{1}^{i_{1}}
...x_{m}^{i_{m}}} } {{i_{1} !...i_{m} !}}}{\frac{{y_{1}^{j_{1}}
...y_{n}^{j_{n}}} } {{j_{1} !...j_{n}
!}}}, \\
 \end{array}
\end{equation}
where $m,n \in {\rm N} \cup {\left\{ {0} \right\}},$ $x$ and $y$
are variables defined in (\ref{eq29}). In the series defined by
(\ref{eq44}) $i_1,...,i_m$ and $j_1,...,j_n$ run from $0$ to
$\infty$.

Evidently, we have

 $$H_{m,m} = F_{{\textbf{A}}}^{(m)} ,
H_{n,0} = F_{{\textbf{B}}}^{(n)} ,H_{2,1} = H_{2} .$$

From the hypergeometric function (\ref{eq44}) we shall define the
following confluent hypergeometric function

$$ \begin{array}{l}{\rm H}_{A}^{(m,n)} \left(
{a,b_1,...,b_m;c_1,...,c_m;x,y} \right) \\= {\mathop {\lim
}\limits_{\varepsilon \to 0}} H_{m + n,m} \left( {a,b_1,...,b_m,
\frac{1}{\varepsilon},...,\frac{1}{\varepsilon};c_1,...,c_m;x,\varepsilon
^{2}y} \right).\\\end{array}
$$

At the determination of the confluent hypergeometric function
${\rm H}_{A}^{(m,n)} $ the equality \cite [p.124]{A} ${\mathop
{\lim }\limits_{\varepsilon \to 0}} \left( {1 / \varepsilon}
\right)_{k} \cdot \varepsilon ^{k} = 1$ ($k$ is a natural number)
has been used. The found confluent hypergeometric function has the
following form

\begin{equation}
\label{eq45} \begin{array}{l}H_{A}^{(m,n)} \left(
{a,b_1,...,b_m;c_1,...,c_m;x,y} \right) \\= {\sum} {\frac{{\left(
{a} \right)_{i_{1} + ... + i_{m} - j_{1} - ... - j_{n}}  \left(
{b_{1}}  \right)_{i_{1}}  ...\left( {b_{m}} \right)_{i_{m}}} }
{{\left( {c_{1}}  \right)_{i_{1}} ...\left( {c_{m}}
\right)_{i_{m}}} } }{\frac{{x_{1}^{i_{1}} ...x_{m}^{i_{m}}} }
{{i_{1} !...i_{m} !}}}{\frac{{y_{1}^{j_{1}} ...y_{n}^{j_{n}}} }
{{j_{1} !...j_{n} !}}},({\left| {x_{1}}  \right|} + ... +
\,{\left| {x_{m}}  \right|} < 1), \\\end{array}
\end{equation}

\bigskip

Now we apply the symbolic operators $\tilde {\nabla} _{x,y}^{m,n}
\left( {h} \right)$ \,and $\tilde {\Delta} _{x,y}^{m,n} \left( {h}
\right)$ to the studying of properties of confluent hypergeometric
function ${\rm H}_{A}^{(m,n)} $ defined by (\ref{eq45}).

Using the formulas (\ref{eq25}) and (\ref{eq27}), we obtain

\begin{equation}
\label{eq46} {\rm H}_{A}^{(m,n)} \left(
{a,b_1,...,b_m;c_1,...,c_m;x,y} \right) = \tilde {\nabla
}_{x,y}^{m,n} \left( {a} \right)F_{A}^{(m)} \left(
{a,b_1,...,b_m;c_1,...,c_m;x} \right)J_{ - a}^{(n)} \left( {y}
\right),
\end{equation}

\begin{equation}
\label{eq47} F_{A}^{(m)} \left( {a,b_1,...,b_m;c_1,...,c_m;x}
\right)J_{ - a}^{(n)} \left( {y} \right) = \tilde {\Delta}
_{x,y}^{m,n} \left( {a} \right){\rm H}_{A}^{(m,n)} \left(
{a,b_1,...,b_m;c_1,...,c_m;x,y} \right),
\end{equation}
where  $J_{a}^{(n)}$ is the Bessel function in $n$ variables:
\begin{equation}
\label{eq43} J_{a}^{(n)} \left( {y} \right): = {\sum\limits_{j_{1}
,...,j_{n} = 0}^{\infty}  {{\frac{{\left( { - 1} \right)^{j_{1} +
...j_{n}}} }{{\left( {1 + a} \right)_{j_{1} + ...j_{n}}} }
}{\frac{{y_{1}^{j_{1}}  ...y_{n}^{j_{n}} }}{{j_{1} !...j_{n}
!}}}}}.
\end{equation}

Now by virtue of formulas (\ref{eq213}) and (\ref{eq214}) from the
formulas (\ref{eq46}) and (\ref{eq47}) we have

\begin{equation}
\label{eq455}
\begin{array}{l}
 {\rm H}_{A}^{(m,n)} \left( {a,b_1,...,b_m;c_1,...,c_m;x,y} \right) = F_{A}^{(m)} \left(
{a,b_1,...,b_m;c_1,...,c_m;x} \right)J_{ - a}^{(n)} \left( {y} \right) \\
 + {\sum\limits_{k = 1}^{\infty}  {{\sum\limits_{l = 1}^{k}
{{\sum\limits_M {{\frac{{\left( { - 1} \right)^{k +
l}k!(k-1)!\left( {b} \right)_{i}}} {{(l-1)!(k-l)!\left( {1 - a}
\right)_{l} \left( {c} \right)_{i}}} }}}} }}
}\frac{x_1^{i_1}...x_m^{i_m}}{i_1!...i_m!}\frac{y_1^{j_1}...y_n^{j_n}}{j_1!...j_n!}
 \\\cdot F_{A}^{(m)} \left( {a +
k,b_1+i_1,...,b_m+i_m;c_1+i_1,...,c_m+i_m;x} \right)J_{ - a +
l}^{(n)} \left( {y} \right), \\
 \end{array}
\end{equation}

\begin{equation}
\label{eq456}
\begin{array}{l}
 F_{A}^{(m)} \left( {a,b_1,...,b_m,;c_1,...,c_m;x} \right)J_{ - a}^{(n)} \left( {y} \right) =
{\rm H}_{A}^{(m,n)} \left( {a,b_1,...,b_m;c_1,...,c_m;x,y} \right) \\
 + {\sum\limits_{k = 1}^{\infty}  {{\sum\limits_{l = 1}^{k}
{{\sum\limits_N }} } {\frac{{k!(k-1)!\left( {a} \right)_{l - k}
\left( {b} \right)_{i}}} {{(l-1)!(k-l)!\left( {1 - a} \right)_{k}
\left( {c} \right)_{i}}} }}
}\frac{x_1^{i_1}...x_m^{i_m}}{i_1!...i_m!}\frac{y_1^{j_1}...y_n^{j_n}}{j_1!...j_n!}
 \\\cdot{\rm H}_{A}^{(m,n)} \left( {a + l -
k,b_1,...,b_m;c_1,...,c_m;x} \right)J_{ - a + k}^{(n)} \left( {y}
\right), \\
 \end{array}
\end{equation}
where
$$ M:=\left\{\left(i_1,...,i_m;j_1,...,j_n\right): i_1\geq0,...,i_m\geq0,j_1\geq0,...,j_n\geq0, i_1+...+i_m=k,j_1+...+j_n=l\right\},$$
 $$ N:=\left\{\left(j_1,...,j_n\right): i_1\geq0,...,i_m\geq0,j_1\geq0, j_1\geq0,...,j_n\geq0,i_1+...+i_m=l, j_1+...+j_n=k\right\}.$$

Thus, we have obtained the decomposition formulas for the multiple
confluent function defined by (\ref{eq45}). We recall that the
multiple Lauricella function $F_{A}^{(m)} $ has an expansion
formula (\ref{e28}).

It is easy to see that in the case when $m=n=1$ the decomposition
formulas (\ref{eq455}) and (\ref{eq456})coincide with the formulas
defined by (\ref{eq330}) and (\ref{eq331}), respectively.

\end{document}